\documentclass[12pt,reqno]{amsart}

\usepackage{xspace}

\usepackage{multicol}

\usepackage[shortlabels]{enumitem}
\setlist[enumerate]{leftmargin=*,font=\upshape,align=parleft,label=(\alph*)}
\setlist[itemize]{leftmargin=*,align=parleft,}

\usepackage{hyperref}
\hypersetup{pdfstartview={XYZ null null 1.00}, pdfpagemode=UseNone, colorlinks,breaklinks,linkcolor=blue,urlcolor=blue, anchorcolor=blue,citecolor=blue}

\usepackage[capitalise]{cleveref}

\usepackage{calrsfs}
\DeclareSymbolFont{defaultmathcal}{OMS}{zplm}{m}{n}
\DeclareSymbolFontAlphabet{\mathcal}{defaultmathcal}
\DeclareSymbolFont{handwritten}{OMS}{rsfs}{m}{n}
\DeclareSymbolFontAlphabet{\handcal}{handwritten}

\usepackage{thmtools}
\usepackage{amssymb}
\usepackage{mathtools}

\usepackage{mdframed}
\newmdenv[
leftmargin = 1cm,
rightmargin = 0pt,
skipabove = 8pt,
skipbelow = 3pt,
innerleftmargin = 8pt,
innertopmargin = 0pt,
innerbottommargin = 0pt,
innerrightmargin = 0pt,
linewidth = 3pt,
topline = false,
rightline = false,
bottomline = false
]{leftbar}

\usepackage{kpfonts}

\usepackage[top=1in, bottom=1in, right=1in, left=1in, bindingoffset=0cm]{geometry}

\usepackage[spacing=true,kerning=true,babel=true]{microtype}

\usepackage[foot]{amsaddr}

\theoremstyle{plain}

\crefname{prop}{Proposition}{Propositions}

\newtheorem{theorem}[equation]{Theorem}

\crefname{cor}{Corollary}{Corollaries}

\theoremstyle{definition}

\crefname{defn}{Definition}{Definitions}
\newtheorem{defn}[equation]{Definition}

\crefname{example}{Example}{Examples}

\crefname{notation}{Notation}{Notations}
\newtheorem{notation}[equation]{Notation}

\theoremstyle{remark}

\crefname{remark}{Remark}{Remarks}
\newtheorem{remark}[equation]{Remark}

\crefname{condition}{Condition}{Conditions}
\newtheorem{condition}[equation]{Condition}

\crefname{claim}{Claim}{Claims}
\newtheorem{claim}[equation]{Claim}

\crefname{caution}{Caution}{Cautions}

\newenvironment{case}[2]{\smallskip\par\noindent \textit{Case}~#1: \rmfamily #2.}{}
\newenvironment{case*}[1]{\smallskip\par\noindent \textit{Case}~#1:\rmfamily}{}

\declaretheoremstyle[
spaceabove=\topsep, 
spacebelow=6pt,
headfont=\normalfont\itshape,
notefont=\normalfont, notebraces={(}{)},
bodyfont=\normalfont,
postheadspace=4pt,
qed=\mbox{\smaller[4]$\boxtimes$}
]{claimproofstyle}
\declaretheorem[name={Proof of Claim}, style=claimproofstyle, unnumbered]{pf}

\newcommand{\0}{\mathbb{\emptyset}}

\newcommand{\N}{\mathbb{N}}

\newcommand{\sR}{{}^\ast\mathbb{R}}

\newcommand{\Pow}{\handcal{P}}

\newcommand{\set}[1]{\left\{ #1 \right\}}

\newcommand*{\defeq}{\mathrel{\vcenter{\baselineskip0.5ex \lineskiplimit0pt \hbox{\scriptsize.}\hbox{\scriptsize.}}}=}

\newcommand{\de}{\delta}
\newcommand{\De}{\Delta}
\newcommand{\e}{\varepsilon}
\renewcommand{\phi}{\varphi}
\newcommand{\si}{\sigma}

\newcommand{\FC}{\mathcal{F}}
\newcommand{\IC}{\mathcal{I}}

\newcommand{\dF}{\mathfrak{d}}

\newcommand{\dom}{\mathrm{dom}}
\newcommand{\im}{\mathrm{im}}

\newcommand{\logX}{\log_{|X|}}

\newcommand{\capsize}[2]{\left|#2\right|_{#1}}

\newcommand{\ldeg}{\de_X}

\newcommand{\Prints}[1][k]{\FC_\pi^{#1}(X)}
\newcommand{\Print}{\boldsymbol{F}}

\newcommand{\Erdos}{Erd\H{o}s\xspace}
\newcommand{\Szemeredi}{Szemer\'{e}di\xspace}

\newcommand{\eqcomment}[2][6pt]{\Big[\text{#2}\Big] \hspace{#1}}

\newcommand{\prel}{\searrow}
\newcommand{\crel}{\nearrow}

\newcommand{\prelG}[1]{\mathrel{\;\:\mathclap{\searrow}{{}^{#1}}\;}}
\newcommand{\crelG}[1]{\mathrel{{}^{#1}\mathclap{\nearrow}{\;\:}}}

\title{A short nonalgorithmic proof of the containers theorem for hypergraphs}

\author{Anton Bernshteyn}
\address[Anton Bernshteyn]{Department of Mathematics, University of Illinois at Urbana-Champaign, IL, 
USA}
\email{bernsht2@illinois.edu}

\author{Michelle Delcourt}
\address[Michelle Delcourt]{School of Mathematics, University of Birmingham, Birmingham, UK
}
\email{m.delcourt@bham.ac.uk}

\author{Henry Towsner}
\address[Henry Towsner]{Department of Mathematics, University of Pennsylvania, PA, 
USA}
\email{htowsner@math.upenn.edu}

\author{Anush Tserunyan}
\address[Anush Tserunyan]{Department of Mathematics, University of Illinois at Urbana-Champaign, IL, 
USA}
\email{anush@illinois.edu}

\thanks{The second author's research was partially supported by EPSRC grant EP/P009913/1 and NSF Graduate Research Fellowship DGE 1144245. The third author's research was partially supported by NSF Grant DMS-1600263. The fourth author's research was partially supported by NSF Grant DMS-1501036.}

\date{}

\begin{document}

\maketitle

\begin{abstract}
	Recently the breakthrough method of hypergraph containers, developed independently by Balogh, Morris, and Samotij~\cite{bms} as well as Saxton and Thomason~\cite{st}, has been used to study sparse random analogs of a variety of classical problems from combinatorics and number theory. The previously known proofs of the containers theorem use the so-called \emph{scythe algorithm}---an iterative procedure that runs through the vertices of the hypergraph. (Saxton and Thomason~\cite{ST_random} have also proposed an alternative, randomized construction in the case of simple hypergraphs.) Here we present the first known deterministic proof of the containers theorem that is not algorithmic, i.e., it does not involve an iterative process. Our proof is less than 4 pages long while being entirely self-contained and conceptually transparent. Although our proof is completely elementary, it was inspired by considering hypergraphs in the setting of nonstandard analysis, where there is a notion of dimension capturing the logarithmic rate of growth of finite sets. Before presenting the proof in full detail, we include a one-page informal outline that refers to this notion of dimension and summarizes the essence of the argument.
\end{abstract}

\section{Introduction}

\noindent\textbf{Hypergraph containers theorems.} An important and extremely active line of research in recent years, especially in combinatorics and number theory, is extending classical results to the so-called ``sparse random setting.'' One breakthrough tool for obtaining such results is the method of hypergraph containers developed independently by Balogh, Morris, and Samotij~\cite{bms} as well as by Saxton and Thomason~\cite{st}. For a high-level overview the hypergraph containers method, see the survey \cite{BMS_survey}.

The hypergraph containers theorem gives a tool for analyzing the structure of all the independent subsets in a hypergraph by ``capturing'' each independent set in one of a small number of ``containers.'' Let $H$ be a $k$-uniform hypergraph with $n$ vertices and $n^{1 + (k-1)\delta}$ edges. In general, $H$ can have close to $2^n$ independent sets, for instance, when all the edges of $H$ span only a small portion of its vertex set. To avoid this, one considers \emph{homogeneous} hypergraphs, i.e., those in which the degree of a vertex cannot significantly exceed the average value $n^{(k-1)\delta}$, and similar upper bounds hold for the codegrees of the sets of $\ell < k$ vertices; for details, see~\cref{defn:hypergraphs:advanced}. According to the containers theorem, if $H$ is sufficiently homogeneous, then each independent set contains a \emph{fingerprint}, which is a subset of size roughly $n^{1 - \de}$. Furthermore, each fingerprint $F$ determines a \emph{container} $C(F)$ of size less than $(1-\alpha)n$, where $\alpha$ is a positive constant, with the property that if $F$ is a fingerprint of an independent set $I$, then $I \subseteq C(F)$. Each container can host at most $2^{(1-\alpha)n}$ independent sets, and the number of containers is bounded by the number of fingerprints, which is at most $2^{o(n)}$, so the total number of independent sets must be much smaller than $2^n$. Note that in this calculation we still used the trivial upper bound $2^{(1-\alpha)n}$ on the number of independent sets inside a given container; in practice, the above approach is usually iterated, leading to particularly strong results.

The containers method has been used to prove (or reprove) sparse versions of theorems originally established for dense hypergraphs. For example, {\Szemeredi's theorem}~\cite{sz} in number theory states that for every $k \in \N$, the largest subset of $[n] \defeq \set{1, \hdots, n}$ which does not contain a $k$-term arithmetic progression (\emph{$k$-AP-free}) is very small, having $o(n)$ elements. This classical result implies that the number of $k$-AP-free subsets of $[n]$ is also small, namely at most $2^{o(n)}$. Considering the hypergraph with vertex set $[n]$ whose edges are the $k$-term arithmetic progressions, the containers method leads to a new combinatorial proof~\cite{bms, st} of a stronger statement known as the {random sparse  version of \Szemeredi's theorem}, which was originally obtained by Schacht~\cite{s} and, independently, by Conlon and Gowers~\cite{cg}. A particularly concise way to phrase this result \cite[Theorem~1.1]{bms} is that for all $m \gg n^{1-1/(k-1)}$, the number of $m$-element $k$-AP-free subsets of $[n]$ is at most ${o(n) \choose m}$. Numerous other applications of the containers method can be found in the original papers \cite{bms, st} as well as in the surveys~\cite{BMS_survey, conlon:icm, RodlSchacht}.

The statements and proofs of the core version of the containers theorem originally appeared in \cite[Proposition 3.1]{bms} and \cite[Theorem 3.4]{st}. We state it here as Theorem~\ref{containers_theorem:finite}. Our main result is a new proof of this theorem, whose advantages are described below. 

\vspace{0.5em}

\noindent\textbf{Our proof.} All previously known proofs of the containers theorem are based on the so-called \emph{scythe algorithm}---in other words, they use an iterative process that runs through the vertices of the hypergraph. (In~\cite{ST_random}, Saxton and Thomason developed a different, \emph{randomized} approach for \emph{simple} hypergraphs, i.e., those in which every pair of vertices lies in at most one edge.) 
In contrast to that, our proof is not algorithmic and provides a deterministic way of building containers in a single step (or, rather, $k$ steps, since it still involves induction on $k$, the uniformity of the hypergraph). It is also conceptually transparent and rather short---under 4 pages. 

Our proof was inspired by an attempt to reprove the containers theorem in the setting of nonstandard analysis, i.e., for ultraproducts of finite hypergraphs. Of course, the theorem for ultraproducts follows from that for finite hypergraphs via the transfer principle, but the present authors were hoping to find a direct proof that would take advantage of the notion of dimension available in the ultraproduct that captures the logarithmic rate of growth. However, it turned out that our approach in the nonstandard setting translated into an even more concise proof for finite hypergraphs, to which we devote the current paper (abandoning ultraproducts altogether).

\vspace{0.5em}

\noindent\textbf{Organization.} The rest of this paper is organized as follows. \cref{sec:terms} establishes standard hypergraph notation and terminology. \cref{sec:statement} begins with our definitions of a \emph{homogeneous hypergraph} and a \emph{print/container pair}, and ends with the statement of the containers theorem in these terms, namely \cref{containers_theorem:finite}. In \cref{sec:proofidea}, we sketch the idea behind our proof inspired by nonstandard analysis. Finally, our proof of \cref{containers_theorem:finite} is presented in \cref{sec:proof}.

\section{Basic notation and terminology}\label{sec:terms}

\noindent The set $\N$ of natural numbers includes $0$ and we denote $\N^+ \defeq \N \setminus \set{0}$. For a set $X$ and $k \in \N^+$, we call a $k$-element subset of $X$ a \emph{$k$-edge} and denote the set of all $k$-edges by $[X]^k$. We refer to a subset $H \subseteq [X]^k$ as a \emph{$k$-uniform hypergraph} (on $X$).

A set $I \subseteq X$ is said to be \emph{$H$-independent} if $H \cap [I]^k = \0$ and we denote
\[
\IC_X(H) \defeq \text{the set of all $H$-independent subsets of $X$}.
\]

\begin{remark}\label{indep_for_1-hypergraphs}
	Note that $[X]^1$ is the set of all singletons of $X$, so for any $1$-hypergraph $H \subseteq [X]^1$, $H$-independent subsets of $X$ are precisely those that are disjoint from $\bigcup H$.
\end{remark}

\begin{notation}\label{notation:hypergraphs:basic}
	Let $X$ be a finite set, $k \in \N^+$, $H \subseteq [X]^k$, and $\ell \in \set{1, \hdots, k-1}$.
	
	\begin{itemize}
		\item For $U \subseteq [X]^\ell$, $V \subseteq [X]^{k - \ell}$, we denote
		\begin{align*}
		[U,V]_H &\defeq \set{e \in H : e = u \cup v \text{ for some } u \in U \text{ and } v \in V},
		\\
		H_U &\defeq \set{v \in [X]^{k-\ell} : u \cup v \in H \text{ for some } u \in U},
		\end{align*}
		We refer to $H_U$ as the \emph{fiber of $H$ over $U$} and if $U = \set{u}$, we write $H_u$ instead of $H_U$. (Another common term for $H_u$ is the \emph{link graph of $u$}.)
		
		\item For each $u \in [X]^\ell$, we denote $\deg_H(u) \defeq |H_u|$.
		
		\item We put $\De_\ell(H) \defeq \max_{u \in [X]^\ell} \deg_H(u)$.
	\end{itemize}
\end{notation}

\begin{notation}
	For sets $A,B$ and a relation $R \subseteq A \times B$, we denote
	\begin{align*}
	\dom(R) &\defeq \set{a \in A : \exists b \in B \text{ with } a R b}
	\\
	\im(R) &\defeq \set{b \in B : \exists a \in A \text{ with } a R b}
	\end{align*}
	and refer to these sets, respectively, as the \emph{domain} and the \emph{image} of $R$.
	
\end{notation}

\section{Statement of the containers theorem}\label{sec:statement}

\noindent Throughout, let $X$ denote a finite nonempty set and $k \in \N^+$.

\begin{defn}\label{defn:hypergraphs:advanced}
	Let $H \subseteq [X]^k$ and $\de \in [0,1]$.
	
	\begin{enumerate}[label=\textup{(\ref*{defn:hypergraphs:advanced}.\alph*)}, itemsep=4pt,leftmargin=*,series=defn:hypergraphs:advanced]
		\item \label{defn:hypergraphs:log-degree} We define the \emph{logarithmic degree} of $H$ as
		\[
		\ldeg(H) \defeq \max_{1 \le \ell < k} \, \frac{1}{k - \ell} \cdot \logX \De_\ell(H).
		\]
		In other words, $\ldeg(H)$ is the least $\de \in [0,1]$ such that $\De_\ell(H) \le |X|^{(k - \ell) \de}$ for all $\ell \in \set{1, \hdots, k-1}$.
		
		\item \label{defn:hypergraphs:bounded} We say that $H$ is \emph{$\de$-bounded} if $\ldeg(H) \le \de$. 
		
		\item \label{defn:hypergraphs:capsize} We let $\capsize{\de}{H}$ denote the maximum size of a $\de$-bounded subhypergraph of $H$, i.e.,
		\[
		\capsize{\de}{H} \defeq \max \set{|H'| : H' \subseteq H\text{ and $H'$ is $\de$-bounded}}.
		\]
		
		\item \label{defn:hypergraphs:homogenous} For $\e > 0$, we say that $H$ is \emph{$(\de, \e)$-homogeneous} if it is $\de$-bounded and $\logX |H| \ge 1 + (k-1) \de - \e$.
	\end{enumerate}
\end{defn}
As the name suggests, $(\delta,\e)$-homogeneity implies that $H$ is ``close to evenly distributed,'' in the sense that for a significant proportion of $x\in X$, $\log_{|X|}|H_x|$ is close to $(k-1)\delta_X(H)$.

\begin{defn}\label{defn:fingerprint}
	Let $\pi \in [0,1]$.
	
	\begin{itemize}
		\item A \emph{$\pi$-fingerprint} (in $X$) is a subset $F \subseteq X$ with $\logX |F| \le \pi$.
		
		\item A \emph{$(\pi, k)$-print} (in $X$) is a sequence $\Print \defeq (F_i)_{0 \le i < \ell}$ of $\pi$-fingerprints (in $X$), where $\ell \leq k-1$. We put $\bigcup \Print \defeq \bigcup_{i < \ell} F_i$ and denote the set of all $(\pi, k)$-prints by $\Prints$.
	\end{itemize}
\end{defn}

\begin{remark}
	In the definition of a $(\pi, k)$-print $\Print$, it is possible that $\ell = 0$ and $\Print = \0$.
\end{remark}

\noindent For $\si \in [0,1]$, we denote $\Pow^\si(X) \defeq \set{C \subseteq X : \logX |X \setminus C| \ge 1 - \si}$.

\begin{defn}\label{defn:print-container}
	Let $k \ge 1$, $H \subseteq [X]^k$, $\pi, \si \in [0,1]$. For relations $\prel \subseteq \IC_X(H) \times \Prints$ and $\crel \subseteq \Prints \times \Pow(X)$, the pair $(\prel,\crel)$ is called a \emph{$(\pi,\si)$-print/container pair} for $H$ if
	\begin{enumerate}[label=\textup{(\ref*{defn:print-container}.\roman*)}, itemsep=4pt,leftmargin=*,series=defn:print-container]
		\item \label{print-container:I_has_P} $\dom(\prel) = \IC_X(H)$;
		
		\item $\dom(\crel) \supseteq \im(\prel)$\label{print-container:P_has_C};
		
		\item \label{print-container:IPC} for each $I \in \IC_X(H)$, $\Print \in \Prints$, and $C \in \Pow(H)$, if $I \prel \Print \crel C$, then 
		\[
		\bigcup \Print \subseteq I \subseteq \bigcup \Print \cup C;
		\]
		
		\item \label{print-container:C_has_cologsize_si} $\im(\crel) \subseteq \Pow^\si(X)$ --- we refer to the sets in $\im(\crel)$ as \emph{containers}.
	\end{enumerate}
\end{defn}

Our main result is a new proof of the following version of the containers theorem:

\begin{theorem}\label{containers_theorem:finite}
	For any $k \in \N^+$, $\pi \in [0,1]$, and $\e > 0$, putting $\de \defeq 1 - \pi$ and $\si \defeq 3^{k-1} \e$, the following holds: For any finite nonempty set $X$ with \[ \e \geq 2k\logX 2 \qquad \text{and} \qquad \pi \geq (k-1)\logX 2,\] any $(\de, \e)$-homogeneous hypergraph $H \subseteq [X]^k$ admits a $(\pi,\si)$-print/container pair.
\end{theorem}

\begin{remark}\label{remark:dependency_on_k}
	We point out that in most applications of the above theorem, $\pi$ and $\de$ are constants independent of $|X|$, while $\e$ and $\si$ are values of order $O(1/\log_2{|X|})$. In particular, saying that for a container $C$, we have $\log_{|X|}|X \setminus C| \geq 1 - \si$, usually means that $|C| \leq (1-\alpha)|X|$ for some positive constant $\alpha$.
\end{remark}

\begin{remark}
	We treat $k$, the uniformity of $H$, as a constant independent of $|X|$. In particular, we make no attempt to optimize the dependence of the ratio $\si/\e$ on $k$. However, some recent applications of the containers method deal with hypergraphs $H$ whose uniformity is a growing function of $|X|$. It would be interesting to see if our approach can lead to improved bounds when $k$ grows with $|X|$, but we do not pursue this question here.
\end{remark}

\section{Idea of proof}\label{sec:proofidea}

\noindent Heuristically, we would like to talk about the ``dimension'' rather than the actual cardinality of the sets appearing in the proof. If the set $X$ has, say, dimension $1$, then the sets whose cardinality has the same ``order of magnitude'' as $|X|$, maybe $|X|/2$ or $|X|/17$, should also have dimension $1$. On the other hand, a set with size $\sqrt{|X|}$ should have dimension $1/2$, while a set with size $|X|^k$ or ${|X| \choose k}$ should have dimension $k$. When $|X|$ is a fixed finite number, this is not well defined. Hence, it makes sense to take a sequence of sets $X_n$ with $|X_n|\rightarrow\infty$, and consider the rates of growth of various sets that appear in the proof. 

This informal idea can be made rigorous by passing to the ultraproduct (as in \cite{arXiv1709.04076, MR1643950}) and working with the \emph{fine pseudofinite dimension} \cite{hrushovski,MR3091666, MR2436141, 2014arXiv1402.5212G}, which captures this property: the dimension of a set is essentially its ``rate of growth relative to $|X_n|$,'' and the dimension is valued in $\sR^+/\mathcal{N}$, where $\sR$ is the ultrapower of the real numbers and $\mathcal{N}$ is the convex subgroup consisting of the ``negligible'' values, namely, those bounded by $\log_{|X|} n$ for some $n \in \N^+$. (Taking the quotient by the negligible values corresponds to identifying the dimension of $Y$ and $Z$ if $|Y| = c|Z|$ for some fixed real number $c$.)

We will now informally outline the proof based on the assumption that a well-behaved notion of dimension exists. We take $\dF_X(H)$ to be the dimension of $|X|^{\delta_X(H)}$; in other words, $\dF_X(H)$ is the coset of $\delta_X(H)$ in the quotient $\sR^+/\mathcal{N}$. Analogous to the definitions above but for a dimension $\dF$, we say $H$ is $\dF$-bounded if, for every $u \in [X]^\ell$,
\[
\dim(H_u) \leq (k - \ell) \dF,
\]
and define $|H|_\dF \defeq \max\set{|H'| : H' \subseteq H \text{ and }H'\text{ is } \dF\text{-bounded}}$.

Then we may attempt to prove our theorem by induction on $k$. The base $k=1$ is easy: just let $I \prel \Print$ when $\Print = \0$ and $\Print \crel C$ when $C = X \setminus \left(\bigcup H\right)$. From now on, assume that $k > 1$. Given $H \subseteq [X]^k$, let $\dF \defeq \dF_X(H)$. For an independent set $I$, we take a maximal fingerprint $F \subseteq I$ such that $H_F$ is \emph{homogeneously expanding}, i.e., 
\[
\dim(|H_F|_{\dF})\geq \dim(F)+(k-1)\dF.
\]
Note that for any $\dF$-bounded $(k-1)$-uniform hypergraph $H'$, we have $\dim(H') \leq 1 + (k-2) \dF$, and hence $\dim(|H_F|_{\dF})\leq 1+(k-2)\dF$, which implies that $\dim(F) \leq 1 - \dF$. If $\dim(F)=1-\dF$, then $H_F$ contains a $\dF$-bounded subhypergraph $G$ with $\dim(G) = 1 + (k-2) \dF$. By the inductive hypothesis, there is a print/container pair $\left(\prelG{\star},\crelG{\star}\right)$ for $G$. Using this pair, we partially define a print/container pair for
$H$ as follows: Put $I \prel (F, F_1, \ldots, F_{\ell-1})$ whenever $I \prelG{\star} (F_1, \ldots, F_{\ell-1})$, and note that $I$, being $G$-independent, must admit such $(F_1, \ldots, F_{\ell-1})$. As for the container relation, we put $(F_0, F_1, \ldots, F_{\ell-1}) \crel C$ whenever $F_0 = F$ and $(F_1, \ldots, F_{\ell-1}) \crelG{\star} C$. Verifying that this works is straightforward.

Now suppose that $\dim(F) < 1 - \dF$. Since $F$ is maximal, for any $x \in I \setminus F$, the fiber $H_{F \cup \set{x}}$ is \emph{not} homogeneously expanding, which, ``morally speaking,'' should mean that
\begin{equation}\label{eq:add_a_point_informal}
	\dim(H_x \setminus H_F) < (k-1) \dF.
\end{equation}
Unfortunately, this is not literally true, and some extra technicalities are necessary to obtain a correct analog of \eqref{eq:add_a_point_informal}, see \S\ref{subsec:container}. Nevertheless, for the purposes of this informal discussion, let us assume that \eqref{eq:add_a_point_informal} actually holds. Then we set $I \prel (F)$ and $(F) \crel C$, where
\[
C \defeq \set{x \in X : \dim(H_x \setminus H_F) < (k-1) \dF_X(H)}.
\]
We certainly have $F \subseteq I \subseteq C \cup F$ and all that remains to check is that $C$ has codimension $1$. We observe that
\[
H \subseteq [X, H_F]_H \cup [C, X^{k-1} \setminus H_F]_H \cup [X \setminus C, X^{k-1}]_H,
\]
and therefore,
\[
\dim(H) \leq \max \set{\dim [X, H_F]_H, \dim [C, X^{k-1} \setminus H_F]_H, \dim [X \setminus C, X^{k-1}]_H}.
\]
But
\begin{align*}
\dim[X, H_F]_H
&\leq 
\dim(H_F) + \dF
\\
&\leq
\dim(F) + (k-1)\dF + \dF
\\
&< 1-\dF + k\dF
\\
&=
1 + (k-1) \dF
=
\dim(H),
\end{align*}
and, using the Fubini property of dimension,
\[
\dim [C, X^{k-1} \setminus H_F]_H \leq \dim(C) + \max_{x \in C} \dim(H_x\setminus H_F) < 1 + (k-1) \dF = \dim(H).
\]
Therefore,
\[
\dim(H) = \dim [X \setminus C, X^{k-1}]_H \leq \dim(X \setminus C) + (k-1) \dF,
\]
which forces $\dim(X \setminus C) = 1$.

For formal reasons, this argument does not quite go through in the rigorous setting of nonstandard analysis: the notion of dimension is ``external'' (not defined by a formula of first-order logic), and therefore such a maximal set $F$ need not exist; in fact it \emph{cannot} exist because adding one point to a nonempty set does not affect its dimension. To fix this, one has to replace the notion of dimension with \emph{logarithmic size}. This is precisely the argument we give below, using bounds on the logarithmic sizes of sets as an approximation to the notion of dimension.

\section{Proof}\label{sec:proof}

\noindent This section is devoted to our proof of \cref{containers_theorem:finite}, so we let $k, \pi, \e, \de, \si, X$ and $H$ be as in its hypothesis and we let $\log$ stand for $\logX$. We adopt the convention that $\log 0 = - \infty$.

We define a $(\pi,\si)$-print/container pair by induction on $k$. For the base case $k = 1$, we let $I \prel \Print$ exactly when $\Print = \0$ and $\Print \crel C$ exactly when $C = X \setminus \left(\bigcup H\right)$. The complement of $C$ is $\bigcup H$ and $\log \left|\bigcup H\right| = \log |H| \ge 1 - \e = 1 - \si$. The rest of the conditions clearly hold as well.

Thus, we may assume that $k > 1$ and that the statement is true for all $1 \le k' < k$. 

\subsection{Choice of constants}\label{subsec:constants}

We take
\setlength{\columnsep}{-2cm}
\begin{multicols}{3}
	\begin{itemize}[labelindent=1em]	
		\item $\de' \defeq \de + \log 2$
		
		\item $\pi' \defeq 1 - \de'$
		
		\item $\tilde\pi \defeq \pi - \e - k \log 2$
		
		\item $\tilde\e \defeq \e + (k+1) \log 2$
		
		\item $\e' \defeq 2\e + 2k\log 2$
		
		\item $\si' \defeq 3^{k-2} \e'$.
	\end{itemize}
\end{multicols}

\noindent Note that since $(k-1)\log 2 \leq \pi$, we have $(k-2)\log 2 \leq \pi'$. Also,
\[
\e' = 2\e + 2k\log 2 \le 
2 \e + 2k \cdot \frac{\e}{2k} \le 2 \e + \e = 3\e,
\]
and hence $\si' \leq \si$.


\begin{defn}
	Call a $\pi$-fingerprint $F$ \emph{expanding} if \[\log|H_F|_{\de'} \geq 1 + (k-2)\de' - \e'.\]
\end{defn}

Notice that a $\pi$-fingerprint $F$ is expanding if and only if the fiber $H_F$ contains a $(\de', \e')$-homogeneous subhypergraph. For each expanding $\pi$-fingerprint $F$, fix an arbitrary $(\de', \e')$-homogeneous subhypergraph $G_F \subseteq H_F$. By the induction hypothesis, $G_F$ admits a $(\pi', \si')$-print/container pair; fix any such $(\pi', \si')$-print/container pair $\big(\prelG{F}, \crelG{F}\big)$.

\subsection{The print relation}

Given $I \in \IC_X(H)$ and $\Print = (F_0, F_1, \hdots, F_{\ell-1}) \in \Prints$, we set $I \prel \Print$ to hold exactly when at least one of the following conditions holds:

\begin{condition}\label{cond:finite:print:inductive}
	We have $\ell \geq 1$, $F_0$ is expanding, $F_0 \subseteq I$, and $I \prelG{F_0} (F_1, F_2, \hdots, F_{\ell-1})$.
\end{condition}

\begin{condition}\label{cond:finite:print:noninductive}
	We have $\ell = 1$, $F_0$ is not expanding, $\log|F_0| < \tilde\pi$, and $F_0$ is maximal among the $\pi$-fingerprints $F$ that are contained in $I$ and satisfy
	\begin{equation}\label{eq:finite:expansive}
	\log\capsize{\de'}{H_F} \ge \log|F| + (k-1) \de' - \tilde\e.
	\end{equation}
\end{condition}

\begin{remark}
	\cref{cond:finite:print:inductive} makes sense, since if $F_0 \subseteq I$ and $G \subseteq H_{F_0}$, then $I$ is $G$-independent.
\end{remark}

\subsection{Condition \labelcref{print-container:I_has_P}}

For a fixed $I \in \IC_X(H)$, there are two cases.

\begin{case}{1}{There is an expanding $\pi$-fingerprint $F \subseteq I$}
	Since $I$ is $G_F$-independent, there is a print $\Print' = (F_1, \hdots, F_{\ell - 1}) \in \Prints[k-1]$ with $I \prelG{F} \Print'$. Therefore, taking $\Print \defeq (F, F_1, \hdots, F_{\ell - 1})$, we see that \cref{cond:finite:print:inductive} holds, so $I \prel \Print$.
\end{case}

\begin{case}{2}{There is no expanding $\pi$-fingerprint $F \subseteq I$}
	
	\begin{claim}\label{claim:there_is_a_fingerprint}
		There is a (possibly empty) set $F \subseteq I$ with $\log|F| < \tilde\pi$ that is maximal among the $\pi$-fingerprints contained in $I$ and satisfying \labelcref{eq:finite:expansive}.
	\end{claim}
	\begin{pf}
		Because $F = \0$ satisfies \labelcref{eq:finite:expansive}, there is a maximal $\pi$-fingerprint $F$ contained in $I$ satisfying \labelcref{eq:finite:expansive}. Then $\log|F| < \tilde\pi$, for otherwise we have \[\log\capsize{\de'}{H_F} \ge \tilde\pi + (k-1) \de' - \tilde\e = (1 - \de' - \e') + (k-1) \de' = 1 + (k-2) \de' - \e',\] which means that $F$ is expanding, contradicting the assumption of our case.
	\end{pf}
	
	\noindent The print $\Print \defeq (F)$, where $F$ is given by \cref{claim:there_is_a_fingerprint}, satisfies \cref{cond:finite:print:noninductive}, so $I \prel \Print$.
\end{case}

\subsection{The container relation}\label{subsec:container}

Given $\Print = (F_0, F_1, \hdots, F_{\ell-1}) \in \Prints$ and $C \in \Pow(X)$, we set $\Print \crel C$ to hold exactly when at least one of Conditions \ref{cond:finite:container:inductive} and \ref{cond:smalldeg} below holds.

\begin{condition}\label{cond:finite:container:inductive}
	We have $\ell \geq 1$, $F_0$ is expanding, and $(F_1, F_2, \hdots, F_{\ell-1}) \crelG{F_0} C$.
\end{condition}

To state \cref{cond:smalldeg}, we need a definition first.

\begin{defn}\label{defn:nabla}
	For $k' \ge 1$, a hypergraph $H' \subseteq [X]^{k'}$, $1 \le t < k'$, and $\de \in [0,1]$, let $\nabla_t^\de(H')$ denote the set of all $u \in [X]^t$ with $\log \deg_{H'}(u) \ge (k'-t) \de$ in $H'$.
\end{defn}

\begin{condition}\label{cond:smalldeg}
	We have $\ell = 1$, $F_0$ is not expanding, and the following holds. Define
	\begin{equation}\label{eq:finite:H-hat}
	H^- \defeq H \setminus \hat{H}, \quad \text{where} \quad \hat{H} \defeq [H_{F_0}, X]_H \cup \bigcup_{t = 1}^{k-2} [\nabla_t^\de(H_{F_0}), [X]^{k-t}]_H.
	\end{equation}
	Then we have
	\begin{equation}\label{eq:finite:C}
	C = \set{x \in X : \log\deg_{H^-}(x) < (k-1) \de' - \tilde\e}.
	\end{equation}
\end{condition}

\subsection{Condition \labelcref{print-container:P_has_C}}

Let $\Print = (F_0, F_1, \hdots, F_{\ell-1}) \in \im(\prel)$. It follows from Conditions \ref{cond:finite:print:inductive} and \ref{cond:finite:print:noninductive} that $\ell \geq 1$.

\begin{case}{1}{$F_0$ is expanding}
	Then \cref{cond:finite:print:inductive} holds. This means that $(F_1, F_2, \hdots, F_{\ell-1}) \in \im(\prelG{F_0}) \subseteq \dom(\crelG{F_0})$. Hence, for some $C \in \Pow(X)$ we have $(F_1, F_2, \hdots, F_{\ell-1}) \crelG{F_0} C$, which yields $\Print \crel C$ by \cref{cond:finite:container:inductive}.
\end{case}

\begin{case}{2}{$F_0$ is not expanding}
	Then \cref{cond:finite:print:noninductive} holds. This means that $\ell = 1$ and there is a (unique) set $C$ satisfying \cref{cond:smalldeg}, so $\Print \crel C$.
\end{case}

\subsection{Condition \labelcref{print-container:IPC}}

We fix $I \in \IC_X(H)$, $\Print = (F_0, F_1, \ldots, F_{\ell-1}) \in \Prints$, and $C \in \Pow(X)$ with $I \prel \Print \crel C$. It follows that $\ell \geq 1$.

\begin{case}{1}{$F_0$ is expanding}
	Set $\Print' \defeq (F_1, F_2, \ldots, F_{\ell-1})$. By the case assumption, Conditions~\ref{cond:finite:print:inductive} and \ref{cond:finite:container:inductive} hold, so $F_0 \subseteq I$ and $I \prelG{F_0} \Print' \crelG{F_0} C$. Therefore, \labelcref{print-container:IPC} applied to $\big(\prelG{F_0}, \crelG{F_0}\big)$ yields $\bigcup \Print' \subseteq I \subseteq \bigcup \Print' \cup C$.
\end{case}

\begin{case}{2}{$F_0$ is not expanding}
	Then Conditions~\ref{cond:finite:print:noninductive} and \ref{cond:smalldeg} hold. In particular, $\ell = 1$. For brevity, let $F \defeq F_0$. By \cref{cond:finite:print:noninductive}, $F \subseteq I$, so it remains to show that each $x \in I \setminus F$ belongs to $C$. Letting $H^-$ be as in \labelcref{eq:finite:H-hat}, we suppose towards a contradiction that $x \notin C$. By \cref{cond:finite:print:noninductive}, $F$ satisfies \labelcref{eq:finite:expansive}, so let $G \subseteq H_F$ be a $\de'$-bounded hypergraph with $\log |G| \ge \log |F| + (k-1) \de' - \tilde\e$.
	
	\begin{claim}
		$G' \defeq G \cup H^-_x$ is $\de'$-bounded.
	\end{claim}
	\begin{pf}
		We fix $\ell \in \set{1, \hdots, k-2}$ and $u \in [X]^\ell$ and show that $\log\deg_{G'}(u) \le (k - 1 - \ell) \de'$. 
		
		If $u \in \nabla_\ell^\de (H_F)$ or $x \in u$, then $G'_u = G_u$, so 
		$
		\deg_{G'}(u) = \deg_G(u) \le |X|^{(k - 1 - \ell) \de'}.
		$
		
		Otherwise, $\deg_{G'}(u) \le \deg_G(u) + \deg_{H^-_x}(u)$. Since $u \notin \nabla_\ell^\de (H_F)$, $\deg_G(u) \le |X|^{(k - 1 - \ell) \de}$. Also, because $x \notin u$,
		\[
		\deg_{H^-_x}(u) = \deg_{H^-} \big(\set{x} \cup u\big) \le \deg_H \big(\set{x} \cup u\big) \le |X|^{(k - 1 - \ell) \de},
		\]
		so $\deg_{G'}(u) \le 2 \cdot |X|^{(k - 1 - \ell) \de} = |X|^{(k - 1 - \ell) \de'}$.
	\end{pf}
	
	Furthermore, $H^-$ and $[H_F, X]_H$ are disjoint, in particular, $H^-_x$ and $H_F \supseteq G$ are disjoint, so
	\begin{align*}
	|G'| &= |G| + |H^-_x| 
	\\
	\eqcomment{Because $x \notin C$}
	&\ge |F| \cdot |X|^{(k-1) \de' - \tilde\e} + |X|^{(k-1) \de' - \tilde\e} 
	\\
	&=
	(|F| + 1) \cdot |X|^{(k-1) \de' - \tilde\e}.
	\end{align*}
	Therefore,
	\[
	\log\capsize{\de'}{H_{F \cup \set{x}}} 
	\ge 
	\log |F \cup \set{x}| + (k-1) \de' - \tilde\e,
	\] 
	i.e., $F \cup \set{x}$ satisfies \labelcref{eq:finite:expansive}. Since $|X|^{\tilde{\pi}} + 1 \leq |X|^\pi$, the set $F \cup \set{x}$ is a $\pi$-fingerprint contained in $I$. This contradicts the properties of $F$ given by \cref{cond:finite:print:noninductive}.
\end{case}

\subsection{Condition \labelcref{print-container:C_has_cologsize_si}}

For a given $C \in \im(\crel)$, fix any $\Print = (F_0, F_1, \hdots, F_{\ell-1})  \in \Prints$ with $\Print \crel C$. It follows that $\ell \geq 1$.

\begin{case}{1}{$F_0$ is expanding}
	Then \cref{cond:finite:container:inductive} holds, so $(F_1, \hdots, F_{\ell-1}) \crelG{F_0} C$, and thus $|X \setminus C| \ge 1 - \si' \ge 1 - \si$.
\end{case}

\begin{case}{2}{$F_0$ is not expanding}
	Then \cref{cond:smalldeg} holds, so $\ell = 1$ and $C$ is defined as in \labelcref{eq:finite:C}. For brevity, let $F \defeq F_0$.
	
	\begin{claim}
		$\log\big|[H_F, X]_H\big| \le \log|F| + k \de$.
	\end{claim}
	\begin{pf}
		$
		\log\big|[H_F, X]_H\big| 
		\le 
		\log|H_F| + \de
		\le
		\log|F| + (k-1) \de + \de
		= 
		\log|F| + k \de.
		$
	\end{pf}
	
	\begin{claim}
		For each $\ell \in \set{1, \hdots, k-2}$, $\log\big|[\nabla_\ell^{\de}(H_F), [X]^{k-\ell}]_H\big| \le \log{k-1 \choose \ell} + \log|F| + k \de$.
	\end{claim}
	\begin{pf}
		Because each edge $e \in H_F$ is counted in the degrees of at most ${k - 1 \choose \ell}$-many points in $[X]^\ell$, we have that
		\[
		\log|\nabla_\ell^{\de}(H_F)| + (k - 1 - \ell) \de \le \log \sum_{u \in [X]^\ell} \deg_{H_F}(u) \le \log {k - 1 \choose \ell} + \log |H_F|.
		\]
		But $\log |H_F| \le \log|F| + (k - 1) \de$, so
		\[
		\log |\nabla_\ell^{\de}(H_F)| 
		\le 
		\log{k-1 \choose \ell} + \log |F| + (k-1) \de - (k - 1 - \ell) \de 
		= 
		\log{k-1 \choose \ell} + \log |F| + \ell \de.
		\]
		Thus,
		\[
		\log\big|[\nabla_\ell^{\de}(H_F), [X]^{k-\ell}]_H\big|
		\le 
		\log \big|\nabla_\ell^{\de} (H_F)\big| + (k - \ell) \de
		\le 
		\log{k-1 \choose \ell} + \log |F| + k \de.
		\qedhere
		\]
	\end{pf}
	
	Let $H^-$ and $\hat{H}$ be defined as in \labelcref{eq:finite:H-hat}.
	\begin{claim}\label{claim:finite:H^-_is_large}
		$\log|H^-| \ge 1 + (k-1) \de - \e - \log 2$.
	\end{claim}
	\begin{pf}
		It follows from the last two claims that
		\begin{align*}
		|\hat{H}| 
		&\le 
		|F| \cdot |X|^{k \de} + \sum_{\ell = 1}^{k-2} {k-1 \choose \ell} \cdot |F| \cdot |X|^{k \de}
		\\
		&=
		\sum_{\ell = 1}^{k-1} {k-1 \choose \ell} \cdot |F| \cdot |X|^{k \de} 
		\\
		&< 
		2^{k-1} \cdot |F| \cdot |X|^{k \de},
		\\
		&<
		2^{k-1} \cdot |X|^{\tilde\pi} \cdot |X|^{k \de}
		\\
		&=
		2^{k-1} \cdot |X|^{\pi-\e-k\log 2} \cdot |X|^{k \de}
		\\
		&=
		2^{k-1} \cdot 2^{-k} \cdot |X|^{1 + (k-1) \de - \e}
		\leq
		\frac{1}{2} \cdot |H|,
		\end{align*}
		so $|H^-| = |H| - |\hat{H}| \ge \frac{1}{2} \cdot |H|$.
	\end{pf}
	
	On the other hand,
	\begin{align*}
	\log \big|[C, [X]^{k-1}]_{H^-}\big| 
	&< 
	|C| + (k-1) \de' - \tilde\e 
	\\
	&= 
	|C| + (k-1) \de + (k-1) \log 2 - \e - (k + 1) \log 2 
	\\
	&\leq 1 + (k-1) \de - \e - 2 \log 2,
	\end{align*}
	and
	\[
	\log \big|[X \setminus C, [X]^{k-1}]_{H^-}\big| \le \log \big|[X \setminus C, [X]^{k-1}]_{H}\big| \le \log |X \setminus C| + (k-1) \de,
	\]
	so, 
	\begin{align*}
	\frac{1}{4} \cdot |X|^{1 + (k-1) \de - \e} + |X \setminus C| \cdot |X|^{(k-1) \de}
	&\ge 
	\big|[C, [X]^{k-1}]_{H^-}\big| + \big|[X \setminus C, [X]^{k-1}]_{H^-}\big|
	\\
	&\ge |H^-| 
	\ge 
	\frac{1}{2}\cdot |X|^{1 + (k-1) \de - \e}.
	\end{align*}
	Therefore, $|X \setminus C| \ge \frac{1}{4} \cdot |X|^{1 - \e}$, so 
	\begin{align*}
	\log |X \setminus C| 
	&\ge 
	1 - \e - 2 \log 2 
	\\
	&\ge 
	1 - \e - 2 \cdot \frac{\e}{2k}
	\\ 
	\eqcomment{Because $k \ge 2$}
	&> 
	1 - 2 \e > 1 - \si.
	\end{align*}
\end{case}

\medskip

The proof of \cref{containers_theorem:finite} is now complete. \hfill\qed

\subsection*{Acknowledgment} We are very grateful to the anonymous referee for helpful comments.

\end{document}